\documentclass{amsart}
\usepackage{amssymb}
\newcommand{\cC}{{\mathcal C}}

\newcommand{\cZ}{{\mathcal Z}}

\newcommand{\cD}{{\mathcal D}}
\newcommand{\cM}{{\mathcal M}}
\newcommand{\be}{{\bf 1}}

\newcommand{\BZ}{{\mathbb Z}}

\newcommand{\BC}{{\mathbb C}}

\newcommand{\sq}{$\square$}
\newcommand{\uomega}{\underline{\omega}}

\newcommand{\Ve}{\mbox{Vec}}
\newcommand{\Rep}{\mbox{Rep}}
\newcommand{\Fun}{\mbox{Fun}}
\newcommand{\Hom}{\mbox{Hom}}
\newcommand{\Mod}{\mbox{Mod}}

\newcommand{\id}{\mbox{id}}

\title{Module categories over the Drinfeld double of a finite group}
\author{Viktor Ostrik}
\email{ostrik@math.mit.edu}
\address{Department of Mathematics, MIT, 77 Massachusetts Ave., Cambridge,
MA 02139}
\thanks{The author was partially supported by NSF grant DMS-0098830}
\date{February 2002}
\begin{document}
\begin{abstract} We classify the module categories over the double (possibly 
twisted) of a finite group. 
\end{abstract} 
\maketitle
\section{Introduction}
Let $\cC$ be a (semisimple abelian) monoidal category. A module category over 
$\cC$ is a (semisimple abelian) category $\cM$ together with a functor 
$\cC \times \cM \to \cM$ and an associativity constraint (= natural 
isomorphism of two composition functors $\cC \times \cC \times \cM \to 
\cM$) satifying natural axioms, see \cite{O}. In physics, one is interested
in the case when $\cC$ is a category
of representations of some vertex algebra $V$ and irreducible objects of
$\cM$ are interpreted as a boundary conditions for the Conformal Field Theory
associated to $V$, see \cite{BPPZ, FS}. Thus it is interesting from a physical
point of view for a given monoidal category $\cC$ to classify all possible 
module categories over $\cC$ (it is known that in many interesting cases the 
list of answers is finite, see \cite{O}). This problem is also of mathematical
interest, for example the module categories with just one isomorphism class of
irreducible objects are exactly the same as the fiber functors $\cC \to \Ve$,
see \cite{O}. It is known that for a fusion categories of $\widehat{sl}(2)$
at positive integer levels the module categories are classified by ADE
Dynkin diagrams, see \cite{Oc, KO, O}. In this note we consider another class
of examples, known in physics as holomorphic orbifold models, see e.g. 
\cite{DVVV, K}. 

Let $G$ be a finite group. It is well known that the monoidal structures 
(= associativity constraints) on the category $\Ve^G$ of $G-$graded vector 
spaces with the usual functor of tensor product are classified by the group
$H^3(G,\BC^*)$, see \cite{MS}. For $\omega \in H^3(G,\BC^*)$ let 
$\Ve^G_\omega$ denote the corresponding monoidal category. Let $\cD(G,\omega)$
be the Drinfeld center of the category $\Ve^G_\omega$, see e.g. \cite{BK, M2}.
The main result of this paper is a classification of module categories over 
the category $\cD(G,\omega)$. The closely related problem of classification of
modular invariants for this category was considered in \cite{CGR}. On the 
other hand the problem of classification of twists (= fiber functors) in the 
Hopf algebra $\cD (G)$ corresponding to $\cD(G,1)$ is interesting from the 
Hopf-algebraic point of view (here 1 is a trivial 3-cocycle).

The paper is organized as follows: in the Section 2 we recall necessary
definitions and basic facts, in Section 3 we state and prove the Main Theorem,
and in Section 4 we consider the simplest nonabelian example --- the case of
the symmetric group in three letters.

I was interested in the problem solved in this note from the beginning of
my work on module categories \cite{O}. But I got an idea that this problem
can be solved (and actually is trivial) only after reading papers \cite{M1,
M2}. So I am grateful to Michael Mueger for this influence. My sincere
gratitude is due to Pavel Etingof for very interesting discussions. I am
greatly indebted to Terry Gannon for providing me with a list of modular 
invariants for the group $S_3$. 

\section{Preliminaries}
In this paper we will consider only abelian semisimple categories over the
field of complex numbers; we assume that all $\Hom-$spaces are finite
dimensional. The rank of such a category is just the number of isomorphism
classes of simple objects (= the rank of the Grothendieck group). We will
consider only categories of finite rank. All functors are assumed to be 
additive. We will consider on many occasions the cohomology of finite
groups with coefficients in $\BC^*$ (with trivial action). We will denote
by 1 the trivial cohomology class. Sometimes we will identify cohomology
class with a representing cocycle; we will always consider only normalized
cocycles (i.e. the value of cocycle is 1 when one of arguments equals to 1).

\subsection{Module categories} Let $\cC =(\cC, \otimes, a, \be, l, r)$ be a 
monoidal category. Here $\cC$ is a category, $\otimes :\cC \times \cC \to \cC$
is a bifunctor, $a: \otimes \circ (\id \times \otimes)\to \otimes \circ 
(\otimes \times \id)$ is a functorial isomorphism of two possible composition 
functors $\cC \times \cC \times \cC \to \cC$, $\be \in \cC$ is the unit 
object, $l: \otimes \circ (\be \times \id)\to \id$ and $r: \otimes \circ (\id 
\times \be)\to \id$ are functorial isomorphisms subject to usual axioms, see 
e.g. \cite{BK}. We will consider only fusion categories, that is rigid 
categories with irreducible unit object, see \cite{ENO}. A (left) module 
category $\cM$ over $\cC$ is a quadruple $(\cM,  \otimes, m, l)$ where 
$\cM$ is a category, $\otimes : \cC \times \cM \to
\cM$ is a bifunctor, $m: \otimes \circ (\otimes \times \id)\to \otimes \circ 
(\id \times \otimes)$ is a functorial isomorphism of two possible composition 
functors $\cC \times \cC \times \cM \to \cM$ and $l: \otimes \circ (\be \times
\id)\to \id$ is a functorial isomorphism of two functors $\cM \to \cM$; these 
data are subject to the pentagon and triangle axioms, see \cite{O}. There are 
obvious notions of a direct sum of module categories, of module functors 
between module categories, see \cite{O}. A module category is called 
indecomposable if it is not module equivalent to a direct sum of nonzero 
categories. One way to construct module categories over $\cC$ is the 
following: let $A$ be an associative algebra in the category $\cC$, see e.g. 
\cite{O}. Then the category $\Mod_\cC (A)$ of right $A-$modules has an obvious
structure of a module category over $\cC$, see {\em loc. cit.} Theorem 3.1 in
{\em loc. cit.} states that any module category over $\cC$ can be obtained in 
this way.

{\bf Example 2.1.} The indecomposable module categories over $\Ve^G_\omega$ 
are classified by the conjugacy classes of pairs $(H, \psi)$ where $H\subset
G$ is a subgroup such that $\omega|_H=1$ and $\psi \in H^2(H,\BC^*)$, see
e.g. \cite{FS, O}. To be more precise let $\uomega \in Z^3(G,\BC^*)$ be a
3-cocycle representing $\omega$. For any $H\subset G$ let $\psi \in 
C^2(H,\BC^*)$ be a 2-cochain such that $d\psi =\uomega|_H$. Then the twisted
group algebra $A(H,\psi)=\BC_\psi[H]$ is an associative algebra in the
category $\Ve^G_\omega$; for any 1-cochain $\chi \in C^1(H,\BC^*)$ the
algebras $A(H,\psi)$ and $A(H,\psi +d\chi)$ are isomorphic. Any semisimple 
indecomposable algebra $A$ in $\Ve^G_\omega$ with $\dim \Hom(\be,A)=1$ is of 
the form $A(H,\psi)$ and thus the set of isomorphism classes of such algebras
is a fibration over $\{ H\subset G | \omega|_H=1\}$ with fiber being a torsor
over $H^2(H,\BC^*)$. In what follows we will often assume that this torsor
is trivialized. 

\subsection{The dual category} Let $\cC$ be a monoidal category and $\cM$ be
an indecomposable module category over $\cC$. Then one defines the dual
category $\cC^*=\Fun_\cC (\cM,\cM)$ as a category of module functors from 
$\cM$ to itself. Note that despite of our notation the category $\cC^*$ 
depends on both $\cC$ and $\cM$. 
In general, semisimplicity of $\cC$ and $\cM$ does not imply that $\cC^*$
is semisimple. But this is always true over a field of characteristic 0
(this follows from the main result of \cite{ENO} in conjunction with 
Theorem 4.2 in \cite{O} and Corollary 6.7 in \cite{N}). Moreover over a field
of characteristic 0 for two module categories $\cM_1$ and $\cM_2$ the
category $\Fun_\cC (\cM_1, \cM_2)$ of module functors from $\cM_1$ to $\cM_2$
is semisimple, see \cite{ENO}. If the category $\cM$ is of the form 
$\Mod_\cC (A)$ for some algebra $A\in \cC$ then one identifies the category 
$\cC^*$ with the category of $A-$bimodules. The operation $\cC \mapsto
\cC^*$ is really a duality, by definition $\cM$ is a module category over
$\cC^*$ and $\cC^{**}=\cC$ (this is an immediate consequence of Theorem 4.2 in
\cite{O}).

{\bf Example 2.2.} Let $\cC =\Ve^G_1$ and $\cM =\Ve$ be the module category
corresponding to the usual forgetful functor $\Ve^G_1\to \Ve$ (this functor
is defined for any $\omega$ but it has a structure of tensor functor if and 
only if $\omega =1$). Then it is easy to see using Theorem 4.2 of \cite{O} 
that $\cC^*=\Rep(G)$ --- the category of representations of the finite group 
$G$. A little bit more generally, let $\cM$ be the module category attached 
to a pair $(H,1)$ in Example 2.1. Then it is easy to see that the category 
$\cC^*$ is equivalent to the category of $H-$equivariant sheaves on $G/H$ 
with tensor product given by the convolution of sheaves.

The following observation is due to M.~Mueger, see \cite{M1}: 

{\bf Proposition 2.1.} The module categories over $\cC$ and over $\cC^*$ are 
in the canonical bijection. 

{\bf Proof.} Let $\cM$ be the module category which connects $\cC$ and $\cC^*$
(that is $\cC^*=\Fun_\cC (\cM,\cM)$). Let $\cM_1$ be a module category over
$\cC$. Then the category $\Fun_\cC (\cM,\cM_1)$ of module functors from $\cM$
to $\cM_1$ is a module category over $\cC^*=\Fun(\cM,\cM)$ via the composition
of functors. Conversely if $\cM_2$ is a module category over $\cC^*$ then
$\Fun_{\cC^*}(\cM,\cM_2)$ is a module category over $\Fun_{\cC^*}(\cM,\cM)=
\cC$. It is obvious that these constructions are inverse to each other. \sq

Note that the category $\Fun(\cM,\cM_1)$ has the following ``elementary''
description: if $\cM=\Mod_\cC(A)$ and $\cM_1=\Mod_\cC(B)$ for some algebras
$A, B \in \cC$ then the category $\Fun(\cM,\cM_1)$ is equivalent to the
category of $(A,B)-$bimodules.  

{\bf Example 2.3.} One instance of this bijection is contained in \cite{AEGN}.
In {\em loc. cit.} the authors define non-commutative ``lowering of indices''.
This can be described as follows: let $\cC$ be a monoidal category (semisimple,
of finite rank) and let $\cM$ be a module category over $\cC$ of rank 1 (=
a fiber functor); thus $\cC$ is the category of representations of a Hopf
algebra and $\cM$ is related with forgetful functor $\cC \to \Ve$; let $\cC^*$
be the dual category of $\cC$ with respect to $\cM$ (thus $\cC^*$ is the
category of representations of the Hopf algebra dual to the previous one). 
One says that a module category $\cM_1$ over $\cC$ of rank 1 (= a fiber 
functor = a twist up to gauge transformation) is {\em non-degenerate} if
the category $\Fun_\cC (\cM, \cM_1)$ is semisimple of rank 1. Theorem 5.10
of \cite{AEGN} says that there is a bijection between non-degenerate module
categories over $\cC$ and over $\cC^*$. From our point of view this is
almost tautological, the bijection being given by an assignment $\cM_1\mapsto
\Fun_\cC (\cM,\cM_1)$. The conditions that $\cM_1$ is of rank 1 and that
$\cM_1$ is non-degenerate are permuted by this bijection. Note that
the assumption of semisimplicity of $\cC$ is not essential for these
considerations. 

\subsection{The Drinfeld center} Let $\cC$ be a monoidal category, $\cM$ be
a module category over $\cC$ and $\cC^*$ be the corresponding dual category.
Then one can consider external product of monoidal categories $\cC \otimes
\cC^*$ and $\cM$ has an obvious structure of a module category over $\cC 
\otimes \cC^*$. The following result is inspired by \cite{M2}:

{\bf Proposition 2.2.} The dual category to $\cC \otimes \cC^*$ with respect
to the module category $\cM$ is equivalent to the Drinfeld center $\cZ (\cC)$ 
of the category $\cC$.

{\bf Proof.} Each functor from $(\cC \otimes \cC^*)^*$ should commute with
$\cC^*$ so is equivalent to the functor $X\otimes ?$ for some $X\in \cC$
(since $\cC^{**}=\cC$). The functor $X\otimes ?$ should commute with all
functors of the form $Y\otimes ?$ for all $Y\in \cC$, so for each $Y\in \cC$
we have an isomorphism of functors $X\otimes (Y\otimes ?)\to Y\otimes (X\otimes
?)$. Clearly this isomorphism comes from the isomorphism of objects $X\otimes
Y\to Y\otimes X$ (by the Yoneda Lemma). We leave for the reader to check that
under this identification the axioms (= commutative diagrams) from the 
definition of a module functors reduce to the axioms of the objects of Drinfeld
center. \sq

We have the following immediate consequence :

{\bf Corollary 2.1.} The categories $\cZ (\cC)$ and $\cZ (\cC^*)$ are 
canonically equivalent.

{\bf Proof.} The construction of Proposition 2.2 is completely symmetric 
with respect to the roles of $\cC$ and $\cC^*$. \sq

{\bf Remark.} A counterpart of this statement
in the theory of Hopf algebras is well known: the Drinfeld doubles of a Hopf
algebra and of its dual are isomorphic.

For a monoidal category $\cC$ let $\cC^{op}$ denote the same category with
opposite tensor product. The following consequence of Proposition 2.2 is 
essentially contained in \cite{M2}:

{\bf Corollary 2.2.} For any monoidal category $\cC$ the module categories 
over $\cZ (\cC)$ are in canonical bijection with module categories over
$\cC \otimes \cC^{op}$.

{\bf Proof.} One should use Proposition 2.2 with $\cM=\cC$ considered as a
module category over $\cC$. \sq

\section{Main results}
\subsection{General setting} Let $G$ be a finite group and let $\omega \in 
H^3(G,\BC^*)$ be a cohomology class. Consider the category $\cC=\Ve^G_\omega$.
For each subgroup $H\subset G$ such that $\omega|_H=1$ and a cohomology
class $\psi \in H^2(H,\BC^*)$ let $A(H,\psi)$ denote the corresponding
algebra in the category $\Ve^G_\omega$, see Example 2.1 (as a $G-$graded 
space $A(H,\psi)$ is supported on $H\subset G$ and each graded component is
one dimensional, the multiplication is given by the 2-cocycle $\psi$;
in other words $A(H,\psi)$ is the group algebra of $H$ twisted by $\psi$).
Consider the monoidal category $\cC (G,\omega,H,\psi)$ of 
$A(H,\psi)-$bimodules in the category $\Ve^G_\omega$.
Equivalently, the category $\cC (G,\omega,H,\psi)$ can be defined as the
category $\cC^*$ with respect to the module category $\cM$ consisting of
$A(H,\psi)-$modules. 

{\bf Theorem 3.1.} The indecomposable module categories over monoidal category
$\cC(G,\omega,H,\psi)$ are parametrized by the conjugacy classes of pairs
$(H_1,\psi_1)$ where $H_1\subset G$ is a subgroup such that $\omega|_{H_1}=1$
and $\psi_1\in H^2(H_1,\BC^*)$ is a cohomology class. Moreover, the module
category corresponding to a pair $(H_1,\psi_1)$ is just the category $\cM$ of
$(A(H,\psi), A(H_1,\psi_1))-$bimodules with tensor product over
$A(H,\psi)$ as a bifunctor $\cC(G,\omega,H,\psi)\times \cM \to \cM$.

{\bf Proof.} This is immediate from the definitions and Proposition 2.1
and Example 2.1. \sq

Now let $(H_1,\psi_1)$ and $(H_2,\psi_2)$ be two pairs consisting of a 
subgroup $H_i\subset G$ such that $\omega|_{H_i}=1$ and a cohomology class
$\psi_i\in H^2(H_i,\BC^*)$. We are going to calculate the number of simple 
objects in the category of $(A(H_1,\psi_1),A(H_2,\psi_2))-$bimodules (so as
a special case we will get the rank of the category $\cC(G,\omega,H,\psi)$).
Let $g\in G$ be an element of $G$. The group $H^g:=H_1\cap gH_2g^{-1}$ has a 
well defined 2-cocycle 
$$\psi^g (h, h'):=\psi_1(h, h')\psi_2(g^{-1}{h'}^{-1}g, g^{-1}h^{-1}g)
\omega(hh'g, g^{-1}{h'}^{-1}g, g^{-1}h^{-1}g)^{-1}\times$$ 
$$\times \omega (h, h', g)\omega(h, h'g, g^{-1}{h'}^{-1}g)$$
(we will see that this is indeed well defined in the proof of Proposition 3.1).

{\bf Proposition 3.1.} Let $(H_1,\psi_1)$ and $(H_2,\psi_2)$ be two pairs as
above. Let $\{ g_i\}_{i\in H_1\setminus G/H_2}$ be a set of reprsentatives of
two-sided $(H_1,H_2)-$cosets in $G$. The rank of the category of $(A(H_1,
\psi_1),A(H_2,\psi_2))-$bimodules equals to
$$\sum_{i\in H_1\setminus G/H_2}m(g_i)$$
where $m(g_i)$ is the number of irreducible projective representations of the 
group $H^{g_i}$ with the center acting via 2-cocycle $\psi^{g_i}$.

{\bf Proof.} Let $M=\oplus_{g\in G}M_g$ be a 
$(A(H_1,\psi_1),A(H_2,\psi_2))-$bimodule. The structure of a bimodule gives us
for each $h_1\in H_1, h_2\in H_2$ the isomorphisms $i^1_{h_1,g}: M_{h_1g}\to
M_g$ and $i^2_{g,h_2}: M_{gh_2}\to M_g$. These isomorphisms should satisfy the
following identities:

1) $i^1_{h_1h_1',g}=\omega(h_1, h_1', g)\psi_1(h_1,h_1')i^1_{h_1',g}\circ
i^1_{h_1,h_1'g}$ (this relation says that $M$ is left $A(H_1,\psi_1)-$module).

2) $i^2_{g,h_2'h_2}=\omega (g, h_2', h_2)^{-1}\psi_2(h_2',h_2)i^2_{g,h_2'}
\circ i^2_{gh_2',h_2}$ (this relation says that $M$ is right 
$A(H_2,\psi_2)-$module).

3) $i^1_{h_1,g}\circ i^2_{h_1g,h_2}=\omega (h_1, g, h_2)i^2_{g,h_2}\circ 
i^1_{h_1,gh_2}$ (that is the structure of left module commutes with the
structure of right module).

There are also some other identities (which we omit) related to the action of
the identity. 

It is clear that the module $M$ is a direct sum of submodules supported on
individual double cosets. If we assume that $M$ is supported on one such coset
and $g$ is a representative of this coset, then it is enough to study the
vector space $M_g$; for each pair $(h_1, h_2)\in H_1\times H_2$ such that
$h_1gh_2=g$ we have an operator $j_{h_1}^g:=i^1_{h_1,g}\circ i^2_{h_1g,h_2}: 
M_g\to M_g$ (note that $h_1=gh_2^{-1}g^{-1}$ so $h_1\in H^g$, also $h_2$ is
completely determined by $h_1$); it is not difficult to see from the above that
$$j_h^g\circ j_{h'}^g=\psi^g(h,h')j_{hh'}^g.$$
This relation is exactly the definition of the $\psi^g-$twisted group algebra 
of $H^g$ and the Proposition is proved. \sq

{\bf Remark.} (i) Actually we will prove that there is a natural bijection 
between simple $(A(H_1,\psi_1),A(H_2,\psi_2))-$bimodules and irreducible
projective representations of groups $H^{g_i}$ with the center acting via
2-cocycle $\psi^{g_i}$.

(ii) The formula for 2-cocycle $\psi^g$ is a generalization of
formula 5.15 from \cite{CGR}.

(iii) The paper of S.~Yamagami \cite{Y} contains results closely related to
Proposition 3.1, see \cite{Y} section 5. Also the result closely related with 
Proposition 3.1 is contained in \cite{AEGN}.

It is easy to deduce from Proposition 3.1 the classification of fiber functors
(= module categories of rank 1) for the category $\cC(G,\omega,H,\psi)$.

{\bf Corollary 3.1.} The fiber functors of the category $\cC(G,\omega,H,\psi)$
are classified by pairs $(H_1,\psi_1)$ where $H_1\subset G$ is a subgroup and
$\psi_1 \in H^2(H_1,\BC^*)$ such that

1) The class $\omega|_{H_1}$ is trivial.

2) The number of double cosets $H\setminus G/H_1$ is 1.

3) The class $\psi|_{H\cap H_1}-\psi_1|_{H\cap H_1}$ is nondegenerate, that is
the twisted group algebra of $H\cap H_1$ with respect to the cocycle 
representing $\psi|{H\cap H_1}-\psi_1|_{H\cap H_1}$ is isomorphic to a matrix
algebra.

{\bf Example 3.1.} Assume that $G=H\cdot H_1$ is an exact factorization, i.e. 
any element of $G$ can be uniquely represented as a product $hh_1$ where
$h\in H, h_1\in H_1$. Then the category $\cC(G,1,H,1)$ has a fiber functor
corresponding to a pair $(H_1,1)$. The corresponding Hopf algebra is called
a Kac algebra and is well known in the theory of Hopf algebras.

\subsection{Main Theorem} Recall that $\cC (G, \omega)$ denotes the Drinfeld
center of the category $\Ve^G_\omega$. Consider the cohomology class 
$\tilde \omega=p_1^*\omega -p_2^*\omega \in H^3(G\times G,\BC^*)$ where
$p_i: G\times G\to G, (g_1,g_2)\mapsto g_i$ are projections. Here is the main 
result of this note:

{\bf Theorem 3.2.} The indecomposable module categories over the
category $\cC (G, \omega)$ are parametrized by the conjugacy classes of pairs 
$(H, \psi)$ where $H\subset G\times G$ is a subgroup such that $\tilde 
\omega|_H=1$ and $\psi \in H^2(H,\BC^*)$.

{\bf Proof.} Let $\Delta (G)\subset G\times G$ denote the diagonal subgroup.
By Corollary 2.2 the category $\cD(G,\omega)$ is equivalent to
the category $\cC(G\times G,\tilde \omega, \Delta(G),1)$. So this Theorem is
a special case of Theorem 3.1. \sq

We will denote the module category corresponding to a pair $(H,\psi)$ by
$\cM (H,\psi)$. Let us consider an action of the group $H$ on the set of
elements of $G$ given by the formula $(h_1,h_2)g=h_1gh_2^{-1}$. Let $\{ g_i
\}_{i\in I}$ be a set of representatives of all orbits under this action.
For each $g_i$ consider a subgroup $H^{g_i}=\{ (h_1,h_2)\in H | h_1g_ih_2^{-1}
=g_i\}$ of $H$; note that since $h_2$ is completely determined by the value of
$h_1$ the group $H^{g_i}$ can be considered as a subgroup of $G$ via the first
projection $H\subset G\times G\to G$. The group $H^{g_i}$ has a well defined
2-cocycle
$$\psi^{g_i}(h,h')=\psi^{-1}((h,g_ihg_i^{-1}),(h',g_ih'g_i^{-1}))\times$$
$$\times \frac{\omega(g_i^{-1},{h'}^{-1},h^{-1})\omega(h,h',{h'}^{-1})\omega
(g_i^{-1}h'g_i,g_i^{-1}{h'}^{-1},h^{-1})}{\omega(hh',{h'}^{-1},h^{-1})\omega
(g_i^{-1}hg_i,g_i^{-1}h'g_i,g_i^{-1}{h'}^{-1}h^{-1})}.$$
Let $m(g_i)$ denote the number of irreducible projective representations
of the group $H^{g_i}$ with center acting via 2-cocycle $\psi^{g_i}$.

{\bf Theorem 3.3.} The rank of the category $\cM (H,\psi)$ is $\sum_{i\in I}
m(g_i)$.

{\bf Proof.} This is an immediate consequence of Proposition 3.1. \sq 

{\bf Example 3.2.} Assume that the cohomology classes $\omega$ and $\psi$ are 
both trivial. Then the category $\cM (H,1)$ is equivalent to the category
of $\Delta (G)-$equivariant sheaves on $G\times G/H$.

{\bf Corollary 3.2.} The fiber functors of the category $\cD(G,\omega)$ are
classified by pairs $(H,\psi)$ where $H\subset G\times G$ is a subgroup and
$\psi \in H^2(H,\BC^*)$ such that

1) The class $\tilde \omega|_H$ is trivial.

2) The number of double cosets $\Delta(G)\setminus G\times G/H$ is 1.

3) The class $\psi|_{H\cap \Delta(G)}$ is nondegenerate, that is the twisted
group algebra of $H\cap \Delta(G)$ with respect to the cocycle representing
$\psi|_{H\cap \Delta(G)}$ is isomorphic to a matrix algebra. 

\subsection{Two unsolved problems} Recall that $\cC=\cD(G,\omega)$ is a 
modular tensor category, see e.g. \cite{BK,M2}. On the other hand, for
a modular tensor category $\cC$ one attaches a {\em modular invariant} to
each module category over $\cC$, see \cite{BEK,O}.

{\bf Problem 1.} Calculate explicitly modular invariants attached to module
categories over $\cD(G,\omega)$.

For a braided tensor category $\cC$ it makes sense to speak about commutative
algebras $A$ in this category, see e.g. \cite{KO}. We will say that a module 
category $\cM$ over $\cC$ is of {\em type I} if $\cM$ is module equivalent to 
$\Mod_\cC (A)$ for a commutative algebra $A$ (in this case the algebra $A$ is
uniquely determined by the module category $\cM$) and $\cM$ is of {\em type II}
otherwise (this is just a translation to the language of module categories
of the corresponding notions due to G.~Moore and N.~Seiberg for modular 
invariants).

{\bf Problem 2.} Which module categories over $\cD(G,\omega)$ are of type I?

One can attach a module category of type I to any subgroup $H\subset G$.
Indeed the category $\Rep(G)$ is a braided subcategory of $\cD(G,\omega)$
and the algebra $A$ of functions on $G/H$ is a commutative algebra in
$\Rep(G)$. But even in a simplest case $G=S_3$ there are some other module
categories of type I, see below. It seems plausible that at least for
$\omega=1$ the type I module categories correspond to pairs $(H,1)$ (so 
cocycle $\psi$ is trivial) such that $H$ is conjugate to a subgroup of 
$G\times G$ invariant under the permutation of factors automorphism of $G\times
G$.

\section{The simplest example}
In this section we consider the simplest non-abelian example, the case of
the group $G=S_3$ --- the symmetric group in three letters. It is well known
that in this case $H^3(G,\BC^*)=\BZ/6\BZ$ (see e.g. \cite{CGR}). For each value
of $\omega \in H^3(G,\BC^*)$ the category $\cD(G,\omega)$ has 8 simple
objects and the same fusion rules (but the structures of monoidal categories
differ). Let $\omega_0$ be a generator of
$H^3(G,\BC^*)$. It is easy to see that $\omega_0$ has a nonzero restriction to
either of the subgroups $\BZ/2\BZ$, $\BZ/3\BZ$ of $G=S_3$. In this section
we will use the additive notations for cohomology groups instead of the
multiplicative notations used above.  

\subsection{Enumeration of subgroups of $G\times G$} In the table below we
enumerate conjugacy classes of all subgroups of $G\times G$ and their 
2-cocycles. The following table is self-explaining; $\Delta$ denotes the
diagonal imbedding $G\to G\times G$, the subgroup $K$ consists of all pairs
of permutations with the same parity. The subgroups $L_1, L_2$ are
$L_1=(\BZ/3\BZ \times e)\Delta (\BZ/2\BZ)$ and $L_2=(e\times \BZ/3\BZ)
\Delta (\BZ/2\BZ)$ (these subgroups were erroneously omitted in the published
version of this note; the mistake was pointed out in \cite{EP}).

$$\begin{array}{|c|c|c|c|c|c|c|} \hline \mbox{notation}&H_1&H_2&H_3&H_4&H_5&
H_6\\ \hline
|H|&1&\multicolumn{3}{|c|}{2}&
\multicolumn{2}{|c|}{3}\\  \hline         
H&e\times e&\BZ/2\BZ \times e&e\times \BZ/2\BZ&\Delta(\BZ/2\BZ)
&\BZ/3\BZ \times e&e\times \BZ/3\BZ\\ \hline          
H^2(H,\BC^*)&0&0&0&0&0&0\\  \hline \end{array}$$

$$\begin{array}{|c|c|c|c|c|c|c|}\hline \mbox{notation}&H_7&H_8&H_9&H_{10}&
H_{11}&H_{12}\\ \hline
|H|&3&4&\multicolumn{4}{|c|}{6}\\ \hline
H&\Delta(\BZ/3\BZ)&\BZ/2\BZ \times \BZ/2\BZ&S_3\times e&e\times S_3&
\Delta(S_3)&\BZ/2\BZ \times \BZ/3\BZ\\ \hline
H^2(H,\BC^*)&0&\BZ/2\BZ&0&0&0&0\\ \hline \end{array}$$

$$\begin{array}{|c|c|c|c|c|c|}\hline \mbox{notation}&H_{13}&H_{14}&H_{15}&
H_{16}\\ \hline
|H|&6&9&\multicolumn{2}{|c|}{12}\\ \hline
H&\BZ/3\BZ \times \BZ/2\BZ&\BZ/3\BZ \times \BZ/3\BZ&S_3\times \BZ/2\BZ&
\BZ/2\BZ \times S_3\\ \hline
H^2(H,\BC^*)&0&\BZ/3\BZ &\BZ/2\BZ&\BZ/2\BZ\\ \hline \end{array}$$

$$\begin{array}{|c|c|c|c|c|c|c|}\hline \mbox{notation}&H_{17}&H_{18}&H_{19}&
H_{20}&H_{21}&H_{22}\\ \hline
|H|&\multicolumn{3}{|c|}{18}&36&6&6\\ \hline
H&S_3\times \BZ/3\BZ&\BZ/3\BZ \times S_3&K&S_3\times S_3&
L_1&L_2\\ 
\hline H^2(H,\BC^*)&0&0&\BZ/3\BZ&\BZ/2\BZ&0&0\\ \hline \end{array}$$

\subsection{Module categories} We say that the subgroup $H$ is admissible for
a 3-cocycle $\omega$ if $\tilde \omega|_H$ (see Theorem 3.2) is trivial. It is
clear that any subgroup is admissible for $\omega =0$ and in this case the 
category $\cD(G,\omega)$ has 28 indecomposable module categories (in both
cases when $H^2(H,\BC^*)=\BZ/3\BZ$ two nontrivial 2-cocycles are permuted by
the action of the normalizer of $H$); if $\omega=\pm \omega_0$ then the 
admissible subgroups are $H_1$, $H_4$, $H_7$, $H_{11}$ and we have 4 
indecomposable module categories; if $\omega=\pm 2\omega_0$ then the 
admissible subgroups are $H_1$, $H_2$, $H_3$, $H_4$, $H_7$, $H_8$, $H_{11}$
and we have 8 indecomposable module categories; finally if $\omega=3\omega_0$
then the admissible subgroups are $H_1$, $H_4$, $H_5$, $H_6$, $H_7$, $H_{11}$,
$H_{14}$, $H_{19}, H_{21}, H_{22}$ and we have 12 indecomposable module 
categories. The following table gives the number of simple objects in module 
categories in all cases (note that the rank of module categories does not 
depend on 2 and 3-cocycles involved since the possible groups $H^g=e, 
\BZ/2\BZ, \BZ/3\BZ, S_3$ have trivial second cohomology, so we present results
only for $\omega=0$; in the other cases one just needs to erase some columns).

$$\begin{array}{|c|c|c|c|c|c|c|c|c|c|c|c|} \hline \mbox{subgroup}&H_1&H_2&H_3&
H_4&H_5&H_6&H_7&H_8&H_9&H_{10}&H_{11}\\ \hline        
|\Delta(G)\setminus G\times G/H|&6&3&3&4&2&2&4&2&1&1&3\\ \hline
\mbox{rank of}\; \cM&6&3&3&6&2&2&10&3&1&1&8\\ \hline \end{array}$$

$$\begin{array}{|c|c|c|c|c|c|c|c|c|c|c|c|} \hline \mbox{subgroup}&H_{12}&
H_{13}&H_{14}&H_{15}&H_{16}&H_{17}&H_{18}&H_{19}&H_{20}&H_{21}&H_{22}\\ 
\hline |\Delta(G)\setminus G\times G/H|&1&1&2&1&1&1&1&2&1&2&2\\ \hline
\mbox{rank of}\; \cM&1&1&6&2&2&3&3&6&3&4&4\\ \hline \end{array}$$

So in particular we see that the untwisted double of $S_3$ has 4 fiber functors
and all the twisted doubles have no fiber functors.

\subsection{Modular invariants} In this section we assume that $\omega =0$ and
give a list of modular invariants attached to the module categories above.
The list of modular invariants for this case was worked out in \cite{CGR} 3.3
(unfortunately the list in a published vesrion of \cite{CGR} is incomplete
and I am grateful to the authors of \cite{CGR} for providing me with a complete
list consisting of 48 invariants). Below we are using the notations for
the simple objects of $\cD(S_3,1)$ from \cite{CGR}: if we identify 
$\cD(S_3,1)$ with the category of $S_3-$equivariant sheaves on $S_3$ with
respect to the action by conjugations then indices 0,1,2 mark objects living
on the unit, 0 corresponds to the trivial representation, 1 to the sign 
representation, 2 to the 2-dimensional irreducible representation of $S_3$; 
indices 3,4,5 mark objects living on the conjugacy class of 3-cycles, 3 
corresponds to the trivial representation and 4,5 to two nontrivial one 
dimensional representations of $\BZ/3\BZ$; indices 6,7 mark objects living
in the conjugacy class of 2-cycles, 6 corresponds to the trivial and 7 to
the nontrivial irreducible representations of $\BZ/2\BZ$. One can compute
the modular invariants attached to the module categories above in the 
following way: first compute the
ranks of the module categories $\cM$ (this was done above) and of the dual
categories $\cC^*$ (= the rank of the category of 
$(A(H,\psi),A(H,\psi))-$bimodules; see the tables below, where $H$ denotes the
module category corresponding to a pair $(H,1)$ and $H^{tw}$ denotes the module
category corresponding to a pair $(H,\psi)$ where $\psi$ is a nontrivial
2-cocycle):
$$\begin{array}{|c|c|c|c|c|c|c|c|c|c|c|c|c|} \hline \mbox{subgroup}&H_1&H_2&
H_3&H_4&H_5&H_6&H_7&H_8&H_8^{tw}&H_9&H_{10}&H_{11}\\ \hline        
\mbox{rank of}\; \cC^*&36&18&18&12&36&36&20&9&9&18&18&8\\ \hline 
\end{array}$$

$$\begin{array}{|c|c|c|c|c|c|c|c|c|c|c|} \hline \mbox{subgroup}&H_{12}&
H_{13}&H_{14}&H_{14}^{tw}&H_{15}&H_{15}^{tw}&H_{16}&H_{16}^{tw}&H_{17}&H_{18}
\\ \hline
\mbox{rank of}\; \cC^*&18&18&36&20&9&9&9&9&18&18\\ \hline \end{array}$$

$$\begin{array}{|c|c|c|c|c|c|c|} \hline \mbox{subgroup}&H_{19}&H_{19}^{tw}&
H_{20}&H_{20}^{tw}&H_{21}&H_{22}\\ \hline        
\mbox{rank of}\; \cC^*&12&8&9&9&12&12\\ \hline \end{array}$$

On the other hand the ranks of $\cM$ and of $\cC^*$ can (at least 
conjecturaly) be computed from the modular invariants, see e.g. \cite{O} Claims
5.3 and 5.4. Now comparing the lists of the module categories and of the
modular invariants one gets some information (unfortunately incomplete, we
can not distinguish for example $H_2$ and $H_{17}$ from this information) on
the correpondence between them. We conjecture that the table below (obtained
using the above information plus some heuristics) gives correct correspondence
between module categories and modular invariants.
 
$$\begin{array}{|c|c|c|c|}\hline 
H_1&|\chi_0+\chi_1+2\chi_2|^2&H_{15}&(\chi_0+\chi_3+
\chi_6)(\chi_0+\chi_2+\chi_6)^*\\ \hline
H_2&(\chi_0+\chi_2+\chi_6)(\chi_0+\chi_1+2\chi_2)^*&H_{15}^{tw}&
(\chi_0+\chi_3)(\chi_0+\chi_2)^*+\chi_1\chi_6^*+\\ \cline{1-2}
H_3&(\chi_0+\chi_1+2\chi_2)(\chi_0+\chi_2+\chi_6)^*&&+\chi_3\chi_6^*+
\chi_6\chi_1^*+\chi_6\chi_2^*+|\chi_7|^2\\ \hline
H_4&|\chi_0+\chi_2|^2+|\chi_1+\chi_2|^2+&H_{16}&(\chi_0+\chi_2+
\chi_6)(\chi_0+\chi_3+\chi_6)^*\\ \cline{3-4}
&+|\chi_6|^2+|\chi_7|^2&H_{16}^{tw}&(\chi_0+
\chi_2)(\chi_0+\chi_3)^*+\chi_1\chi_6^*+\\ \cline{1-2}
H_5&(\chi_0+\chi_1+2\chi_3)(\chi_0+\chi_1+2\chi_2)^*&&+\chi_2\chi_6^*+
\chi_6\chi_1^*+\chi_6\chi_3^*+|\chi_7|^2\\ \hline
H_6&(\chi_0+\chi_1+2\chi_2)(\chi_0+\chi_1+2\chi_3)^*&H_{17}&(\chi_0+\chi_3+
\chi_6)(\chi_0+\chi_1+2\chi_3)^*\\ \hline
H_7&|\chi_0+\chi_1|^2+2|\chi_2|^2+2|\chi_3|^2+&H_{18}&
(\chi_0+\chi_1+2\chi_3)(\chi_0+\chi_3+\chi_6)^*\\ \cline{3-4}
&+2|\chi_4|^2+2|\chi_5|^2&H_{19}&
|\chi_0+\chi_3|^2+|\chi_1+\chi_3|^2+\\ \cline{1-2}
H_8&|\chi_0+\chi_2+\chi_6|^2&&+|\chi_6|^2+|\chi_7|^2\\ \hline
H_8^{tw}&|\chi_0+\chi_2|^2+\chi_1\chi_6^*+\chi_2\chi_6^*+&H_{19}^{tw}&
|\chi_0|^2+|\chi_1|^2+\chi_2\chi_3^*+\chi_3\chi_2^*\\
&+\chi_6\chi_1^*+\chi_6\chi_2^*+|\chi_7|^2&&+|\chi_4|^2+|\chi_5|^2+
|\chi_6|^2+|\chi_7|^2\\ \hline
H_9&(\chi_0+\chi_3+\chi_6)(\chi_0+\chi_1+2\chi_2)^*&H_{20}&|\chi_0+\chi_3+
\chi_6|^2\\ \hline
H_{10}&(\chi_0+\chi_1+2\chi_2)(\chi_0+\chi_3+\chi_6)^*&H_{20}^{tw}&
|\chi_0+\chi_3|^2+\chi_1\chi_6^*+\chi_3\chi_6^*+\\ \cline{1-2}
H_{11}&|\chi_0|^2+|\chi_1|^2+|\chi_2|^2+|\chi_3|^2+&&+\chi_6\chi_1^*+
\chi_6\chi_3^*+|\chi_7|^2 \\ \cline{3-4}
&+|\chi_4|^2+|\chi_5|^2+|\chi_6|^2+|\chi_7|^2&H_{21}&(\chi_0+\chi_3)(\chi_0+
\chi_2)^*+|\chi_6|^2+\\ \cline{1-2}
H_{12}&(\chi_0+\chi_2+\chi_6)(\chi_0+\chi_1+2\chi_3)^*&&(\chi_1+\chi_3)
(\chi_1+\chi_2)^*+|\chi_7|^2\\ \hline
H_{13}&(\chi_0+\chi_1+2\chi_3)(\chi_0+\chi_2+\chi_6)^*&H_{22}&(\chi_0+\chi_2)
(\chi_0+\chi_3)^*+|\chi_6|^2+\\ \cline{1-2}
H_{14}&|\chi_0+\chi_1+2\chi_3|^2&&(\chi_1+\chi_2)
(\chi_1+\chi_3)^*+|\chi_7|^2\\ \hline
H_{14}^{tw}&|\chi_0+\chi_1|^2+2\chi_2\chi_3^*+
2\chi_3\chi_2^*+&&\\
&+2|\chi_4|^2+2|\chi_5|^2&&\\ \hline
\end{array}$$

\end{document}